\documentclass[12pt,a4paper]{amsart}
\usepackage{amssymb,amsmath}

\headsep=1cm \numberwithin{equation}{section}
\newtheorem{theorem}{Theorem}[section]
\newtheorem{lemma}[theorem]{Lemma}

\newtheorem{corollary}[theorem]{Corollary}

\theoremstyle{definition}
\newtheorem{definition}[theorem]{Definition}
\theoremstyle{remark}
\newtheorem{remark}[theorem]{Remark}
\newtheorem{example}[theorem]{Example}

\newcommand{\Ass}{\operatorname{Ass}}
\newcommand{\im}{\operatorname{im}}
\newcommand{\grade}{\operatorname{grade}}

\newcommand{\E}{\operatorname{E}}

\newcommand{\V}{\operatorname{V}}
\newcommand{\D}{\operatorname{D}}

\newcommand{\Ext}{\operatorname{Ext}}
\newcommand{\Supp}{\operatorname{Supp}}
\newcommand{\Tor}{\operatorname{Tor}}
\newcommand{\Hom}{\operatorname{Hom}}

\newcommand{\lo}{\longrightarrow}
\newcommand{\fm}{\frak{m}}
\newcommand{\fp}{\frak{p}}

\newcommand{\fa}{\frak{a}}
\newcommand{\fb}{\frak{b}}

\begin{document}
\author[Divaani-Aazar and Mafi ]{Kamran Divaani-Aazar and Amir Mafi}
\title[Associated primes of local cohomology modules of ...]
{Associated primes of local cohomology modules of weakly Laskerian
modules}

\address{K. Divaani-Aazar, Department of Mathematics, Az-Zahra University,
Vanak, Post Code 19834, Tehran, Iran-and-Institute for Studies in
Theoretical Physics and Mathematics, P.O. Box 19395-5746, Tehran,
Iran.} \email{kdivaani@ipm.ir}

\address{A. Mafi,  Institute of Mathematics, University for Teacher
Education, 599 Taleghani Avenue, Tehran 15614, Iran.}



\begin{abstract} The notion of weakly Laskerian modules was
introduced recently by the authors. Let $R$ be a commutative
Noetherian ring with identity, $\fa$ an ideal of $R$, and $M$ a
weakly Laskerian module. It is shown that if $\fa$ is principal,
then the set of associated primes of the local cohomology module
$H_{\fa}^i(M)$ is finite for all $i\geq 0$. We also prove that
when $R$ is local, then $\Ass_R(H_{\fa}^i(M))$ is finite for all
$i\geq 0$ in the following cases: (1) $\dim R\leq 3$, (2) $\dim
R/\fa\leq 1$, (3) $M$ is Cohen-Macaulay and for any ideal $\fb$,
with $l=\grade(\fb,M)$, $\Hom_R(R/\fb,H_{\fb}^{l+1}(M))$ is weakly
Laskerian.
\end{abstract}

\maketitle

\section{Introduction}

Throughout, let $R$ be a commutative Noetherian ring with
identity. Let $\fa$ be an ideal of $R$ and $M$ an $R$-module. For
each $i\geq 0$, the $i$-th local cohomology module of $M$ with
respect to $\fa$ is defined as
$$H_{\fa}^i(M)=\underset{n}{\varinjlim}\Ext_R^i(R/\fa^n,M).$$
The reader can refer to [{\bf 3}] for basic properties of local
cohomology. In [{\bf 7}], Hartshorne defined an $R$-module $M$ to
be $\fa$-cofinite if $\Supp_RM\subseteq \V(\fa)$ and
$\Ext_R^i(R/\fa,M)$ is finitely generated for all $i\geq 0$. He
then asked when the local cohomology modules of a finitely
generated module are $\fa$-cofinite. In this regard, the best
known result is that for a finitely generated $R$-module $M$ if
either $\fa$ is principal or $R$ is local and $\dim R/\fa=1$, then
the $H_{\fa}^i(M)$'s are $\fa$-cofinite. These results are proved
in [{\bf 9}, Theorem 1] and [{\bf 4}, Theorem 1], respectively.

It is easy to see that an $\fa$-cofinite module has only finitely
many associated primes. Huneke [{\bf 8}] raised the following
question: If $M$ is a finitely generated $R$-module, then the set
of associated primes of $H_{\fa}^i(M)$ is finite for all ideals
$\fa$ of $R$ and all $i\geq 0$. Singh [{\bf 16}] gives a
counter-example to this question (also see [{\bf 17}] for some
other counter-examples). On the other hand, Brodmann and Lashgari
[{\bf 2}, Theorem 2.2] showed that if for a finitely generated
$R$-module $M$ and an integer $t$, the local cohomology module
$H_{\fa}^i(M)$ is finitely generated for all $i<t$, then
$\Ass_R(H_{\fa}^t(M))$ is finite. For another proof of this
result, see [{\bf 10}]. As in [{\bf 6}], an $R$-module $M$ is
called weakly Laskerian if $\Ass_R(M/N)$ is finite for all
submodules $N$. Let $M$ be a weakly Laskerian module and $t\in
\mathbb{N}$. In [{\bf 6}], we proved the set of associated primes
of the first non $\fa$-cofinite local cohomology module of $M$ is
finite. That clearly implies the result mentioned above due to
Brodmann and Lashgari. Also regarding the Artinianness of local
cohomology, it is shown in [{\bf 5}], that in many cases weakly
Laskerian modules behave similar to finitely generated modules.

In this article, we continue studying the set of associated primes
of local cohomology of weakly Laskerian modules. Example 2.3 shows
that the class of weakly Laskerian modules is much larger than
that of Noetherian modules. Our main aim in this paper is to show
that in conjunction with finiteness properties of local
cohomology, in many cases weakly Laskerian modules behave similar
to finitely generated modules. By applying the known techniques we
can get the analogues of some finiteness results for finitely
generated modules for weakly Laskerian modules .

In section 2, we first review weakly Laskerian modules. Then we
introduce the notion of weak cofiniteness and show that many
results concerning cofiniteness of finitely generated modules are
valid for weak cofiniteness of weakly Laskerian modules. In
particular, it is proved that the Change of Ring Principle for
weak cofiniteness holds.

The main results of this article appear in section 3. We prove
that the local cohomology modules of a weakly Laskerian module are
weakly cofinite in several cases. As a result, we deduce that the
sets of associated primes of the local cohomology modules of a
weakly Laskerian module are finite in these cases. For instance,
it is shown that if $\fa$ is a principal ideal of $R$ and $M$ a
weakly Laskerian $R$-module, then  $\Ass_R(H_{\fa}^i(M))$ is
finite for all $i\geq0$. Also, we prove that if $R$ is a local
ring and $M$ a weakly Laskerian $R$-module, then
$\Ass_R(H_{\fa}^i(M))$ is finite for all
$i\geq0$, in the following cases:\\
a) $\dim R\leq 3$,\\
b) $\dim R/\fa\leq 1$.\\
Finally, we show that if $M$ is a Cohen-Macaulay module over a
local ring $R$ such that for any ideal $\fb$ with
$l=\grade(\fb,M)$, $\Hom_R(R/\fb,H_{\fb}^{l+1}(M))$ is weakly
Laskerian, then $\Ass_R(H_{\fa}^i(M))$ is finite for all $i\geq0$.
This extends the main result of [{\bf 1}].

\section{Weakly cofinite modules}

In this section, we first recall the definition of weakly
Laskerian modules. Then we bring a lemma which is needed in the
sequel.

\begin{definition} i) An $R$-module $M$ is said to be {\it Laskerian}
if any submodule of $M$ is an intersection of a finite number of
primary submodules.\\
ii) (See [{\bf 6}]) An $R$-module $M$ is said to be {\it weakly
Laskerian} if the set of associated primes of any quotient module
of $M$ is finite.
\end{definition}

Obviously, any Noetherian module is Laskerian and it is clear that
any Laskerian module is weakly Laskerian. So, any Noetherian
module is weakly Laskerian. We need to the following lemma in the
sequel.

\begin{lemma} i) Let  $0\lo L\lo M\lo N\lo 0$ be an exact
sequence of $R$-modules. Then $M$ is weakly Laskerian if and only
if $L$ and $N$ are both weakly Laskerian. Thus any subquotient of
a weakly Laskerian module as well as any finite direct sum of
weakly Laskerian modules is weakly Laskerian.\\
ii) Let $M$ and $N$ be two $R$-modules. If $M$ is weakly Laskerian
and $N$ is finitely generated, then $\Ext^i_R(N,M)$ and
$\Tor^R_i(N,M)$ are weakly Laskerian for all $i\geq 0$.\\
iii) Let $M$ be an $R$-module such that $\Supp_R M$ is finite.
Then $M$ is weakly Laskerian. In particular, any Artinian
$R$-module is weakly Laskerian.
\end{lemma}

{\bf Proof.} The proof of i) is easy and we leave it to the
reader.\\
ii) See [{\bf 6}, Lemma 2.3].\\
iii) Let $N$ be a submodule of $M$. We have
$$\Ass_R(M/N)\subseteq \Supp_R(M/N)\subseteq \Supp_R M.$$ Thus the
set of associated primes of any quotient module of $M$ is finite.
$\Box$

Now, we provide several examples for this class of modules.

\begin{example} i) An $R$-module $M$ is said to be minimax if $M$ has
a finitely generated submodule $S$ such that $M/S$ is Artinian
(see [{\bf 18}]). By Lemma 2.2, it follows that any minimax
$R$-module is weakly Laskerian.

ii) Let $E$ be the minimal injective cogenerator of $R$. If for an
$R$-module $M$ the natural map from $M$ to $\Hom_R(\Hom_R(M,E),E)$
is an isomorphism, then $M$ is said to be Matlis reflexive. Also,
an $R$-module $M$ is said to be linearly compact if each system of
congruences
$$x\equiv x_i (M_i),$$ indexed by a set $I$ and where the $M_i$
are submodules of $M$, has a solution $x$ whenever it has a
solution for every finite subsystem. It is known that if either
$M$ is reflexive or linearly compact, then $M$ has a finitely
generated submodule $S$ such that $M/S$ is Artinian (see e.g.
[{\bf 6}, Example 2.2]). Hence any reflexive $R$-module and also
any linearly compact $R$-module is weakly Laskerian.
\end{example}

Next, we present a generalization of the notion of cofiniteness.

\begin{definition} Let $\fa$ be an ideal of $R$ and $M$ an
$R$-module. We say that $M$ is $\fa$-{\it weakly cofinite} if
$\Supp_R M\subseteq \V(\fa)$ and $\Ext_R^i(R/\fa,M)$ is weakly
Laskerian for all $i\geq 0$.
\end{definition}

\begin{example} i) Let $\fa$ be an ideal of $R$ and $M$ an $R$-module
with $\Supp_R M\subseteq \V(\fa)$. If $M$ is weakly Laskerian,
then by Lemma 2.2 ii), it turns out that $M$ is $\fa$-weakly
cofinite. In particular, if either $M$ is finitely generated or
Artinian, then $M$ is $\fa$-weakly cofinite.

ii) Every $\fa$-cofinite module is $\fa$-weakly cofinite.\\
\end{example}

We collect some important known properties of cofinite modules in
a lemma. Then we study the corresponding properties for weakly
cofinite modules.

\begin{lemma} i) Suppose $M$ is a finitely generated $R$-module with
 $\Supp_R M\subseteq \V(\fa)$. Then $M$ is $\fa$-cofinite.\\
 ii) If $0\lo L\lo M\lo N\lo 0$ is an exact sequence of
 $R$-modules and two of the modules in the exact sequence are
 $\fa$-cofinite, then so is the third one.\\
 iii) If $M$ is a $\fa$-cofinite $R$-module, then $\Ass_RM$ is finite.\\
 iv) For every $\fa$-cofinite $R$-module $M$, the $R$-module $M/\fa M$ is
 finitely generated.
\end{lemma}

{\bf Proof.} i) is clear by the definition of cofiniteness.\\
ii) follows easily from the long Ext sequence induced by the exact
sequence $$0\lo L\lo M\lo N\lo 0,$$ and the fact that if, in an
exact sequence, two of the modules in the exact sequence are
Noetherian,
then the third one is also Noetherian.\\
iii) Since $\Supp_R M\subseteq \V(\fa)$, it follows that
$\Ass_RM=\Ass_R(0:_M\fa).$ As $M$ is $\fa$-cofinite, the module
$(0:_M\fa)\cong \Ext^0_R(R/\fa,M)$ is finitely generated. Thus
$\Ass_RM$ is finite.\\
iv) See [{\bf 13}, Corollary 1.2]. $\Box$

\begin{remark} i) Example 2.5 i) shows that the analogue of
Lemma 2.6 i) holds for weak cofiniteness. Also, in view of Lemma
2.2 i), it follows that the analogue of Lemma 2.6 ii)
holds for weak cofiniteness.\\
ii) Because the set of associated primes of a weakly Laskerian
module is finite, by using the same argument as in the proof of
Lemma 2.6 iii), we can deduce that the set of associated
primes of a $\fa$-weakly cofinite module is finite.\\
iii) Recall that a module $M$ is said to have finite Goldie
dimension if $M$ does not contain an infinite direct sum of
nonzero submodules, or equivalently, the injective envelope
$\E(M)$ of $M$ decomposes as a finite direct sum of indecomposable
injective submodules. By [{\bf 13}, Proposition 1.3], the Goldie
dimension of an $\fa$-cofinite module is finite. This is not the
case for $\fa$-weakly cofinite modules. To this end, take
$M=\oplus_{i\in \mathbb{N}}R/\fm$, where $\fm$ is a fixed maximal
ideal of $R$. Since $\Supp_R M=\{\fm \}$ by Lemma 2.2 iii), it
turns out that $M$ is weakly Laskerian. Therefore $M$ is
$\fm$-weakly cofinite, although its Goldie dimension is not
finite.
\end{remark}

Next, we are going to prove the analogue of Lemma 2.6 iv) for weak
cofiniteness, but we first need some preliminary results. In view
of Lemma 2.2, by using the same proof as that used in [{\bf 4},
Proposition 1], we can deduce the following.

\begin{lemma} Let $\fa$ be an ideal of $R$ and $N$ an $R$-module.
For any given integer $n$, the following are equivalent:\\
i) $\Ext^i_R(R/\fa,N)$ is weakly Laskerian for all $i\leq n$.\\
ii) $\Ext^i_R(M,N)$ is weakly Laskerian for all finitely generated
$R$-module $M$ with $\Supp_R M\subseteq \V(\fa)$, and all $i\leq
n$.
\end{lemma}

The Change of Ring Principle for cofiniteness was proved by
Delfino and Marley (see [{\bf 4}, Proposition 2]). In the sequel,
we prove the Change of Ring Principle for weak cofiniteness. The
proof is an adaption of the proof of [{\bf 4}, Proposition 2].

\begin{theorem} Let the ring $T$ be a homomorphic image of $R$. Let $\fa$
be an ideal of $R$ and $M$ a $T$-module. Then $M$ is  $\fa$-weakly
cofinite as an $R$-module if and only if $M$ is $\fa T$-weakly
cofinite as a $T$-module.
\end{theorem}

{\bf Proof.} Assume $T=R/I$ for some ideal $I$ of $R$ and let $N$
be a $T$-module. Then we have $$\Ass_TN=\{\fp/I: \fp\in
\Ass_RN\}$$ and $$\Supp_TN=\{\fp/I: \fp\in \Supp_RN\}.$$ Thus $N$
is weakly Laskerian as a $T$-module if and only if  it is weakly
Laskerian as an $R$-module. Also, we deduce that $\Supp_R
M\subseteq \V(\fa)$ if and only if $\Supp_T M\subseteq \V(\fa T)$.
By [{\bf 15}, Theorem 11.65], there is a Grothendieck spectral
sequence
$$E_2^{p,q}=\Ext_T^p(\Tor_q^R(T,R/\fa),M)\underset{p}{\Longrightarrow}
\Ext_R^{p+q}(R/\fa,M).$$ First assume that $M$ is $\fa T$-weakly
cofinite. Because $\Supp_T(\Tor_q^R(T,R/\fa))\subseteq \V(\fa T)$,
by Lemma 2.8, it turns out that $E_2^{p,q}$ is a weakly Laskerian
$T$-module for all $p$ and $q$. For each $r\geq 2$, $E_r^{p,q}$ is
a subquotient of $E_2^{p,q}$. Hence, by Lemma 2.2 i) $E_r^{p,q}$
is a weakly Laskerian $T$-module for all $r\geq 2$ and all
$p,q\geq 0$. There is an integer $r\gg 0$ such that
$E_{\infty}^{p,q}=E_r^{p,q}$ for all $p,q\geq 0$. Also, for each
$n\in \mathbb{N}_0$ there is a bounded filtration
$$0=\phi^{n+1}H^n\subseteq \phi^nH^n\subseteq \dots \subseteq
\phi^1H^n \subseteq \phi^0H^n=H^n$$ for the module
$H^n=\Ext_R^n(R/\fa,M)$ such that $E_{\infty}^{p,n-p}\cong
\phi^pH^n/\phi^{p+1}H^n$ for all $p=0,1,\dots ,n.$ Now, by using
Lemma 2.2 i) successively, we deduce that $\Ext_R^n(R/\fa,M)$ is a
weakly Laskerian $R$-module for all $n\geq 0$. That is $M$ is
$\fa$-weakly cofinite.

Conversely, assume that $M$ is $\fa$-weakly cofinite. By induction
on $n$, we show that $E_2^{n,0}=\Ext_T^n(T/\fa T,M)$ is a weakly
Laskerian $T$-module for all $n\geq 0$. Since the module
$$\Hom_T(T/\fa T,M)\cong \Hom_R(R/\fa,M)$$ is weakly Laskerian,
the claim holds for $n=0$. Now, assume that $n> 0$ and that
$E_2^{p,0}$ is weakly Laskerian for all $p< n$. By Lemma 2.8, it
turns out that $E_2^{p,q}$ is weakly Laskerian for all $p< n$  and
$q\geq 0$. We have $E_r^{n,0}\cong E_{\infty}^{n,0}$ for all $r\gg
0$. As $H^n=\Ext_R^n(R/\fa,M)$ is weakly Laskerian, it follows
that for each $p=0,1\dots ,n,$ $E_{\infty}^{p,n-p}$ is weakly
Laskerian. From the induction hypothesis, we get that $\im
d_{r-1}^{n-r+1,r-2}$ is weakly Laskerian. Hence $\ker
d_{r-1}^{n,0}$ is weakly Laskerian. By continuing this argument,
we deduce that $E_2^{n,0}=\ker d_2^{n,0}$ is weakly Laskerian.
$\Box$

Now, we are ready to present the last result of this section. Its
proof is a slight modification of the proof of [{\bf 13},
Corollary 1.2]

\begin{theorem} Let $\fa$ be an ideal of $R$ and $M$ a
$\fa$-weakly cofinite $R$-module. Then $M/\fa M$ is weakly
Laskerian.
\end{theorem}

{\bf Proof.} Let $a_1,a_2,\dots ,a_n$ be generators of $\fa$. Let
$\Theta :R[X_1,X_2,\dots X_n]\lo R$ be the natural ring
homomorphism defined by $\Theta(f)=f(a_1,a_2,\dots ,a_n)$. Then
$\Theta$ is surjective and $(X_1,X_2,\dots X_n)R=\fa$. Denote the
ideal $(X_1,X_2,\dots X_n)$ by $\fb$. By Theorem 2.9,
$\Ext^i_{R[X_1,X_2,\dots X_n]}(R[X_1,X_2,\dots X_n]/\fb,M)$ is
weakly Laskerian $R[X_1,X_2,\dots X_n]$-module for all $i\geq 0$.
Let $K_\bullet(X_1,X_2,\dots X_n)$ denote the Koszul complex of
the ring $R[X_1,X_2,\dots X_n]$ with respect to $X_1,X_2,\dots
X_n$. Since $X_1,X_2,\dots X_n$ is a regular sequence in the ring
$R[X_1,X_2,\dots X_n]$, it follows that $K_{.}(X_1,X_2,\dots X_n)$
is a free resolution for $R[X_1,\dots ,X_n]/\fb$. Thus
$$\Ext^n_{R[X_1,X_2,\dots X_n]}(R[X_1,X_2,\dots X_n]/\fb,M)\cong
M/\fb M.$$ Therefore $M/\fa M$ is a weakly Laskerian $R$-module.
$\Box$

\section{Associated primes of local cohomology modules}

In this section, we prove that the local cohomology modules of a
weakly Laskerian module are weakly cofinite in several cases.

Let $M$ be a finitely generated $R$-module and $t$ a nonnegative
integer. In [{\bf 12}, Proposition 2.5], it is shown that if
$H_{\fa}^i(M)$ is $\fa$-cofinite for all $i\neq t$, then
$H_{\fa}^t(M)$ is also $\fa$-cofinite. In the following we show
that the same result holds for weak cofiniteness. Then, as a
result we deduce several results concerning finiteness of
associated primes of local cohomology modules.

\begin{theorem} Let $\fa$ be an ideal of $R$ and $M$ a weakly
Laskerian module. Suppose there exists an integer $t\geq 0$ such
that $H_{\fa}^i(M)$ is $\fa$-weakly cofinite for all $i\neq t$.
Then $H_{\fa}^t(M)$ is also $\fa$-weakly cofinite.
\end{theorem}

{\bf Proof.} By [{\bf 15}, Theorem 11.38], there is a spectral
sequence
$$E_2^{p,q}:=\Ext_R^p(R/\fa,H^q_{\fa}(M))\underset{p}{\Longrightarrow}
\Ext_R^{p+q}(R/\fa,M).$$ For each $r\geq 2$, we consider the exact
sequence $$0\lo \ker d_r^{p,t}\lo E_r^{p,t}\overset{d_r^{p,t}}\lo
E_r^{p+r,t-r+1}. (*)$$ It follows from the hypotheses that the
$R$-module $E_r^{p+r,t-r+1}$ is weakly Laskerian. Note that
$E_r^{p,q}$ is a subquotient of $E_2^{p,q}$ for all $p,q\in
\mathbb{N}_{0}$.

There is an integer $s$ such that $E_{\infty}^{p,q}=E_r^{p,q}$ for
all $r\geq s$. Also, for each $n\in \mathbb{N}_0$ there is a
bounded filtration
$$0=\phi^{n+1}H^n\subseteq \phi^nH^n\subseteq \dots \subseteq
\phi^1H^n \subseteq \phi^0H^n=H^n$$ for the module
$H^n=\Ext_R^n(R/\fa,M)$ such that $E_{\infty}^{p,n-p}\cong
\phi^pH^n/\phi^{p+1}H^n$ for all $p=0,1,\dots ,n.$ Thus
$E_{\infty}^{p,q}$ is weakly Laskerian for all $p,q$. Since
$$E_s^{p,t}=\ker d_{s-1}^{p,t}/\im d_{s-1}^{p-s+1,t+s-2},$$ it
follows that $\ker d_{s-1}^{p,t}$ is weakly Laskerian. Hence by
using the exact sequence $(*)$ for $r=s-1$, we deduce that
$E_{s-1}^{p,t}$ is weakly Laskerian. By continuing this argument
repeatedly for integers $s-1,s-2,\dots ,3$ instead of $s$, we
obtain that $E_2^{p,t}$ is weakly Laskerian for all $p\geq 0$.
This completes the proof. $\Box$

By using the above result, we can deduce the following corollary.

\begin{corollary} Let $\fa$ be an ideal  of a local ring $(R,\fm)$
and $M$ a finitely generated $R$-module. Assume that $\fa$
contains an $M$-filter regular sequence $x_1,\dots ,x_t$ and that
$H_{\fa}^i(M)$ is $\fa$-weakly cofinite for all $i>t$. Then
$H_{\fa}^t(M)$ is $\fa$-weakly cofinite.
\end{corollary}

{\bf Proof.} By [{\bf 14}, Theorem 3.1], $H_{\fa}^i(M)$ is
Artinian for all $i<t$. Hence, the claim follows by Lemma 2.2 iii)
and Theorem 3.1. $\Box$

The following extends the main result of [{\bf 9}].

\begin{corollary} Let $\fa$ be a principal ideal of $R$ and $M$ a
weakly Laskerian module. Then $H_{\fa}^i(M)$ is $\fa$-weakly
cofinite for all $i\geq 0$.
\end{corollary}

{\bf Proof.} Since $H_{\fa}^0(M)$ is a submodule of $M$, it turns
out that $H_{\fa}^0(M)$ is weakly Laskerian, by Lemma 2.2 i).
Also, $H_{\fa}^i(M)=0$ for all $i>1$. Hence $H_{\fa}^i(M)$ is
$\fa$-weakly cofinite for all $i\neq 1$. Therefore, the claim
follows by Theorem 3.1. $\Box$

\begin{corollary} Let $(R,\fm)$ be a local ring of dimension $d$.
Let $\fa$ be an ideal of $R$.

a) For any $R$-module $M$, the modules $H_{\fa}^d(M)$ and
$H_{\fa}^{d-1}(M)$  are weakly Laskerian.

b) Assume that $M$ is a weakly Laskerian $R$-module. The following
assertions hold.\\
i) if $d\leq 3$, then $H_{\fa}^i(M)$ is $\fa$-weakly cofinite for
all $i\geq 0$.\\
ii) if $d=4$, then $H_{\fa}^1(M)$ is $\fa$-weakly cofinite if and
only if $H_{\fa}^2(M)$ is $\fa$-weakly cofinite.\\
iii) if $d=4$ and the $\fa$-ideal transform module of $M$ is
weakly Laskerian, then  $H_{\fa}^i(M)$ is $\fa$-weakly cofinite
for all $i\geq 0$.
\end{corollary}

{\bf Proof.} a) is immediate, because by [{\bf 11}, Corollaries
2.3 and 2.4] $\Supp_RH_{\fa}^d(M)$ and $\Supp_RH_{\fa}^{d-1}(M)$
are finite.
Now, we prove b).\\
i) and ii) are clear by Theorem 3.1 and a).\\
iii) Let $\D_{\fa}(M)$ denote the $\fa$-ideal transform module of
$M$. There is an exact sequence $$0\lo H_{\fa}^0(M)\lo M\lo
\D_{\fa}(M)\lo H_{\fa}^1(M)\lo 0.$$ Hence by Lemma 2.2 i),
$H_{\fa}^0(M)$ and $H_{\fa}^1(M)$ are weakly Laskerian. On the
other hand, by a) $H_{\fa}^4(M)$ and $H_{\fa}^3(M)$ are both
weakly Laskerian. Therefore the conclusion follows, by Theorem
3.1. $\Box$

\begin{example} i) Let $(R,\fm)$ be a local ring of dimension
$3$ and let $x,y,z$ be a system of parameters of $R$. Let
$\fa=(xz,yz)$. Then, by [{\bf 13}, Theorem 2.2], $H_{\fa}^2(R)$ is
not $\fa$-cofinite. But $H_{\fa}^2(R)$ is $\fa$-weakly cofinite,
by Corollary 3.4 i).

ii) Let $R$ and $\fa$ be either as in Example 3.7 or as in Example
3.10 in [{\bf 12}]. Then $\Hom_R(R/\fa,H_{\fa}^2(R))$ is not
finitely generated. Thus $H_{\fa}^2(R)$ is not $\fa$-cofinite,
while in each case $H_{\fa}^i(R)$ is $\fa$-weakly cofinite for all
$i\geq 0$, by Corollary 3.4 i).

iii) Let $k$ be a field, $R=k[x,y][[u,v]]$, $\fa=(u,v)$, and
$f=ux+vy$. Then $H_{\fa}^2(R/(f))$ is not $\fa$-cofinite (see
[{\bf 7}, \S 3]). However, by Corollary 3.4 i), $H_{\fa
}^i(R/(f))$ is $\fa(R/(f))$-weakly cofinite for all $i\geq 0$. In
particular, it follows that the set of associated primes of
$H_{\fa}^i((R/(f))$ is finite for all $i\geq 0$.

iv) Let $R=\mathbb{Z}[x,y,z,u,v,w]$ and $f=ux+vy+wz$. For the
ideal $\fa=(x,y,z)$, Singh [{\bf 16}] has shown that the set of
associated primes of $H_{\fa}^3(R/(f))$ is not finite. Thus
$H_{\fa}^3(R/(f))$ is not $\fa$-weakly cofinite, by Remark 2.7
ii).
\end{example}

The following extends [{\bf 4}, Theorem 1] in some sense.

\begin{theorem} Let $(R,\fm)$ be a local ring and $\fa$ an ideal
of $R$ with $\dim R/\fa\leq 1$. Let $M$ be an $R$-module. Then
$H_{\fa}^i(M)$ is weakly Laskerian for all $i\geq 0$.
\end{theorem}

{\bf Proof.} Since $\dim R/\fa\leq 1$, it follows that $\V(\fa)$
is contained in the set of minimal prime ideals of $\fa$ union
with $\{\fm\}$. Thus $$\Supp_R(H_{\fa}^i(M))\subseteq
\V(\fa)\subseteq \Ass_RR/\fa\cup \{\fm\},$$ for all $i\geq 0$.
Hence, by Lemma 2.2 iii) $H_{\fa}^i(M)$ is weakly Laskerian.
Therefore, the assertion follows. $\Box$

The following extends [{\bf 1}, Theorem 1.2].

\begin{theorem} Let $(R,\fm)$ be a local ring and $M$ a
Cohen-Macaulay $R$-module. Suppose that for any ideal $\fa$ of
$R$, with $t=\grade(\fa,M)$, $\Hom_R(R/\fa,H_{\fa}^{t+1}(M))$ is
weakly Laskerian. Then $\Ass_R(H_{\fa}^i(M))$ is a finite set for
any ideal $\fa$ of $R$ and all $i\geq 0$.
\end{theorem}

{\bf Proof.} Let $i\geq 0$ be a fixed integer and $\fa$ an ideal
of $R$. We show that the set of associated primes of the local
cohomology module $H_{\fa}^i(M)$ is a finite set. If $i\leq
\grade(\fa,M)$, then the claim follows by [{\bf 2}, Theorem 2.2].
Hence assume that $i>\grade(\fa,M)$. We may and do assume that
$H_{\fa}^i(M)\neq 0$. By [{\bf 1}, Lemma 2.4], there exists an
ideal $\fb\supseteq \fa$ of $R$ such that $\grade(\fb,M)=i-1$ and
$H_{\fb}^i(M)\cong H_{\fa}^i(M)$. Now, by the assumption
$\Hom_R(R/\fb,H_{\fb}^i(M))$ is weakly Laskerian. Thus
$\Ass_R(H_{\fa}^i(M))$ is a finite set, as required. $\Box$

\begin{theorem} Let $(R,\fm)$ be a local ring and $M$ a
Cohen-Macaulay $R$-module. Let $t\in \mathbb{N}_{0}$ be a fixed
integer such that $H_{\fa}^t(M)$ is weakly Laskerian for any ideal
$\fa$ of $R$. Then $\Ass_R(H_{\fa}^{t+1}(M))$ is a finite set for
any ideal $\fa$ of $R$.
\end{theorem}

{\bf Proof.} Let $\fa$ be an ideal of $R$. As in the proof of
Theorem 3.7, we may assume that $t+1\geq \grade(\fa,M)$. Also, we
can and do assume that $H_{\fa}^{t+1}(M)\neq 0$. Since $M$ is a
Cohen-Macaulay $R$-module, it follows by [{\bf 1}, Lemma 2.4],
that there exists an ideal $\fb\supseteq \fa$ of $R$ such that
$\grade(\fb,M)=t$ and $H_{\fb}^{t+1}(M)\cong H_{\fa}^{t+1}(M)$. By
the assumption the module $H_{\fb}^t(M)$ is weakly Laskerian.
Thus, it follow from [{\bf 6}, Theorem 2.5] that
$\Ass_R(H_{\fa}^{t+1}(M))$ is a finite set. $\Box$

%%%%%%%%%%%%%%%%%%%%%%%%%%%%%%%%%%%%%%%%%%%%%%%%%%%%%%%%%%%%%%%%%%%%%%%%%%%%%

\end{document}